\documentclass[12pt]{amsart}
\usepackage[latin1]{inputenc}
\usepackage[english]{babel}
\usepackage{amsfonts,amssymb,amsmath,latexsym,stmaryrd,upref,a4}
\usepackage{graphicx}

\input xypic
\xyoption{all}
\CompileMatrices

\usepackage{color}

\newtheorem{Thm}{Theorem}[section]
 \newtheorem{Rem}[Thm]{Remark}
 
\newtheorem{Ex}[Thm]{Example}
\newtheorem{Prop}[Thm]{Proposition}

\begin{document}

\title
[Indecomposables]{Indecomposable finite-dimensional representations of a class of Lie algebras and Lie superalgebras}

\author{Hans Plesner Jakobsen}
\address{
  Department of Mathematical Sciences\\ University of Copenhagen\\Universitetsparken 5\\
   DK-2100, Copenhagen,
  Denmark} \email{jakobsen@math.ku.dk}\date{\today}

 \maketitle

\section{Introduction}

The topic of indecomposable finite-dimensional representations of the
Poincar\'e group was first studied in a systematic way by S. Paneitz
(\cite{pan1}\cite{pan2}). In these investigations only representations
with one source were considered, though by duality, one representation
with 2 sources was implicitly present.  

The idea of nilpotency was mentioned indirectly in Paneitz's articles,
but a more down-to-earth method was chosen there. 

The results form a part of a major investigation by S. Paneitz and
I. E. Segal into physics based on the conformal group. Induction from
indecomposable representations plays an important part in this
theory. See (\cite{psv}) and references cited therein. 

The defining representation of the Poincar\'e group, when given as a
subgroup of SU(2,2) (see below), is indecomposable. This
representation was studied by the present author prior to the articles
by Paneitz in connection with a study of special aspects of Dirac
operators and positive energy  representations of the conformal group
(\cite{jak}).  

Indecomposable representations in theoretical physics have also been
used in a major way in a study by G. Cassinelli, G. Olivieri,
P. Truini, and V. S. Varadarajan (\cite{raja}). The main object is the
Poincar\'e group. In an appendix to the article, the indecomposable
representations of the 2-dimensional Euclidean group are considered,
and many results are obtained. This group can also be studied by our
method, but we will not pursue this here. One small complication is
that the circle group is abelian.

In the article at hand, we wish to sketch how, by utilizing nilpotency
to its fullest extent while using methods from the theory of universal
enveloping algebras, a complete description of the indecomposable
representations may be reached. In practice, the combinatorics is
still formidable, though. 

It turns out that the method applies to both a class of ordinary Lie
algebras and to a similar class of Lie superalgebras. 

Besides some examples, due to the level of complexity we will only
describe a few precise results. One of these is a complete
classification of which ideals can occur in the enveloping algebra of
the translation subgroup of the Poincar\'e group. Equivalently, this
determines all indecomposable representations with a single,
1-dimensional source. Another result is the construction of an
infinite-dimensional family of inequivalent representations already in
dimension 12. This is much lower than the 24-dimensional
representations which were thought to be the lowest possible. The
complexity increases considerably, though yet in a manageable fashion,
in the supersymmetric setting.  Besides a few examples,  only a
subclass of ideals of the enveloping algebra of the super Poincar\'e
algebra will be determined in the present article.

\section{Finite-dimensional indecomposable representations of the Poincar\'e group}

\label{fd}

We are here only interested in what happens on the level of the Lie
algebra. Equivalently, we consider a double covering of the Poincar\'e
group given by

\begin{displaymath}
  P=\left\{\left(\begin{array}{cc}a&0\\k&{a^\star}^{-1}\end{array}\right)
\mid
    a\in SL(2,{\mathbb R})\; ; k^\star=a^\star k a^{-1} \in
    gl(2,{\mathbb C})\right\}.
\end{displaymath}

This is a subgroup of  $SU(2,2)$ when the latter is defined by the hermitian form

\begin{displaymath}
 \left(\begin{array}{cc}0&i\\-i&0\end{array}\right) .
\end{displaymath}

For our purposes, we may equivalently even consider the group 

\begin{displaymath}
  P=\left\{\left(\begin{array}{cc}u&0\\z&v\end{array}\right)
\mid
    u,v\in SL(2,{\mathbb C})\; ; z \in
    gl(2,{\mathbb C})\right\}.
\end{displaymath}

Let $G_0$ denote the group $SL(2,{\mathbb C})\times SL(2,{\mathbb C})$. Thus,
\begin{displaymath}
  G_0=\left\{\left(\begin{array}{cc}u&0\\0&v\end{array}\right)
\mid
    u,v\in SL(2,{\mathbb C})\;\right\}.
\end{displaymath}
For what we shall be doing, it does not matter if we work with this
group, its Lie algebra, or with $SU(2)\times SU(2)$. 

It is important to consider the abelian  Lie algebra

\begin{displaymath}
  {\mathfrak p}^-=\left\{\left(\begin{array}{cc}0&0\\z&0\end{array}\right)
\mid
     z \in
    gl(2,{\mathbb C})\right\}.
\end{displaymath}
along with its enveloping algebra $${\mathcal U}( {\mathfrak
  p}^-)={\mathcal S}( {\mathfrak p}^-) ={\mathbb C}[z_1,z_2,z_3,z_4].$$ 
The last equality comes from writing the $2\times 2$ matrix $z$
above as
$$z=\left(\begin{array}{cc}z_1&z_2\\z_3&z_4\end{array}\right).$$

In passing we make the important observation that the polynomial $\det z=z_1z_4-z_2z_3$ is
invariant in the sense that $\det uzv=\det z$ for all $u,v\in
SL(2,{\mathbb C})$.

\bigskip

We make the basic assumption that all representations and equivalences
are over ${\mathbb C}$. This has the powerful consequence that the
abelian algebra  ${\mathfrak p}^-$ acts nilpotently. 

The general setting is the following: We consider a reductive Lie
algebra ${\mathfrak g}_0$ in which the elements of the abelian ideals
are given by semi-simple elements and a nilpotent Lie algebra
${\mathfrak n}$ together with a Lie algebra homomorphism $\alpha$ of
${\mathfrak g}_0$ into the derivations $Der({\mathfrak n})$ of
${\mathfrak n}$. This gives rise, in the usual fashion, to the
semi-direct product
\begin{displaymath}
	{\mathfrak g}={\mathfrak g}_0\times_\alpha {\mathfrak n}.
\end{displaymath}

In this situation a well-known result  from algebra
(\cite{jac},\cite{hu}) can easily be generalized to include the
${\mathfrak g}_0$ invariance. 

\medskip 
Recall that a flag in a vector space $V$ is a sequence of subspaces
${0}=W_0\subsetneq W_1\subsetneq \dots\subsetneq W_r=V$.

\begin{Thm}\label{lem1} Suppose given a representation of ${\mathfrak
    g}$ in some finite-dimensional vector space $V$. Then there is a
  flag of subspaces such that ${\mathfrak n}$ maps  $W_i$ into
  $W_{i-1}$ for $i=1,\dots,r$ and such that each $W_i$ is invariant
  and completely decomposable under ${\mathfrak g}_0$. 
\end{Thm}

\medskip

We associate a graph to the indecomposable representation $V$ as
follows: The nodes are the ${\mathfrak g}_0$ irreducible
representations that occur. Two nodes, labeled by irreducibles $V_{1}$
and $V_2$, are connected by an arrow pointing from $V_1$ to $V_2$ if
$V_2\subset{\mathfrak n}^-\cdot V_{1}$ inside $V$. If there are
multiplicities in the isotypic components the situation becomes more
complicated. If the multiplicity at the node $i$ is $n_i$ one can
simply place $n_i$ black dots at the node. They can be placed in a
stack or in a circle. In case $n_i>1$, $n_j>1$ there may also be a
number $a_{i,j}>1$ of arrows from $i$ to $j$, and this needs also to
be indicated. The simplest way is just to attach the $n_i$ to each
node and to attach the $a_{i,j}$ to the arrow from $i$ to $j$ with the
further stipulation that only numbers strictly greater that 1 need to
be given.  We shall not pursue such details here; see, however, the
third of the simple examples below.

The theorem above has the immediate consequence that there are no closed paths in this graph.

\begin{Rem}The assumption of finite dimension is essential
  here. Already on the level of the polynomial algebra ${\mathbb
    C}[z_1,z_2,z_3,z_4]$, if one quotients out by the ideal generated
  by an inhomogeneous polynomial in  $\det z$ there will be closed
  loops, but the resulting module is infinite-dimensional. If one
  insists on finite-dimensionality, one must have all homogeneous
  polynomials in the quotient after a certain degree. Thus $\det z^n$
  for some $n$ must be in the ideal. This precludes an inhomogeneous
  polynomial in $\det z$ since in ${\mathbb C}[T]$ (where  $T=\det
  z$), any inhomogeneous polynomial $p(T)$ is relatively prime to
  $T^n$ for any $n=0,1,2,\dots$. 
\end{Rem}

Given any such directed graph, any node with arrows only pointing out
is as usual called a source. The opposite is called a sink. There is a
simple way whereby one may reverse all arrows, thereby turning sources
into sinks, and vice versa: Replace $V$ by its dual module
$V^\prime$. 

\medskip

Simple situations:

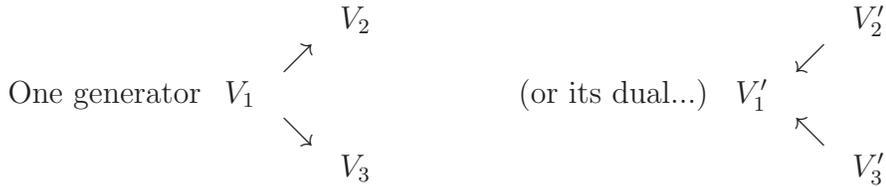
\begin{displaymath}
  \textrm{ One generator } \begin{array}{rcl}&&V_2\\&\nearrow
    \\V_1\\ &\searrow \\&&V_3\end{array}\qquad\qquad\textrm{
    (or its dual...) } \begin{array}{rcl}&&V_2^\prime\\&\swarrow
    \\V_1^\prime\\ &\nwarrow \\&&V_3^\prime\end{array}\qquad\qquad 
\end{displaymath}

This leads to decomposable representations if the targets (sinks) are equal
(respectively if the origins (sources) are equivalent). Otherwise they are indecomposable.
\bigskip
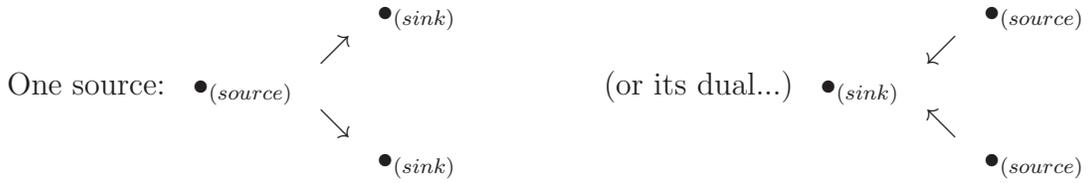
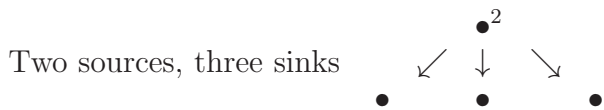
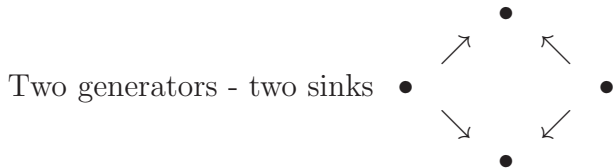
\begin{displaymath}
  \textrm{ One source: } \begin{array}{rcl}&&\bullet_{(sink)}\\&\nearrow
    \\\bullet_{(source)}\\ &\searrow \\&&\bullet_{(sink)}\end{array}\qquad\qquad\textrm{
    (or its dual...) }
  \begin{array}{rcl}&&\bullet_{(source)}\\&\swarrow
    \\\bullet_{(sink)}\\ &\nwarrow
    \\&&\bullet_{(source)}\end{array}\qquad\qquad 
\end{displaymath}

\begin{displaymath}
  \textrm{Two sources, three sinks } \ \begin{array}{rcllc}&&\bullet^2\\&\swarrow
    &\downarrow&\searrow\\\bullet&&\bullet&&\bullet\end{array}\qquad\qquad\qquad\qquad\qquad
\end{displaymath}
\bigskip

\begin{displaymath}
  \textrm{ Two generators - two sinks } \begin{array}{rcllc}&&\bullet\\&\nearrow
    &&\nwarrow\\\bullet&&&&\bullet\\ &\searrow &&\swarrow
    \\&&\bullet\end{array}\qquad\qquad\qquad\qquad\qquad 
\end{displaymath}

\bigskip

\subsection{One generator} We consider only the Poincar\'e algebra. Let $V_0$ denote an irreducible
finite-dimensional representation of ${\mathfrak g}_0=sl(2,{\mathbb
  C})\times sl(2,{\mathbb C})$, given by non-negative integers $(n,m)$
so that the spins are $(\frac{n}2,\frac{m}2)$ and the dimension is $(n+1)(m+1)$. 

Let $\Pi$ denote an indecomposable finite-dimensional representation
of  ${\mathcal S}( {\mathfrak p}^-)\times_s {\mathfrak g}_0$ in a space $V_\Pi$, generated by
a ${\mathfrak g}_0$ invariant source $V_0$. Let ${\mathcal P}(V_0)$ denote the space of polynomials in the
variables $z_1,z_2,z_3,z_4$ with values in $V_0$. This is generated by
polynomials of the form $p_0(z_1,z_2,z_3,z_4)\cdot v$ for $v\in V_0$
and $p_0$ a complex polynomial. We consider this a left ${\mathcal
  S}({\mathfrak p}^-)\times_s {\mathfrak g}_0$  module in the obvious
way. The map  
 
\begin{equation}
  p_0(z_1,z_2,z_3,z_4)\cdot v\mapsto \pi(p_0)v
\end{equation}
is clearly ${\mathcal S}( {\mathfrak p}^-)\times_s {\mathfrak g}_0$ equivariant.

The decomposition of the ${\mathfrak g}_0$ module  ${\mathcal S}( {\mathfrak
  p}^-)$ is well-known and is given by the representations $\textrm{spin}(\frac{n}2,\frac{n}2)$ for
$n=0,1,2,\dots$. Each occurs with infinite multiplicity due to the
invariance of $\det z$ under ${\mathfrak g}_0$. We will describe these representations in detail below.

\medskip

The decomposition of  ${\mathcal P}(V_0)$ into irreducible ${\mathfrak g}_0$
representations follows easily from this using the well-known
decomposition of the tensor product of two irreducible representations
of $su(2)$. The decomposition of  ${\mathcal P}(V_0)$ is in general
more degenerate than what results from the invariance of $\det z$.

It is clear that there exists a sub-module ${\mathcal I}\subseteq {\mathcal P}(V_0)$ such that
\begin{displaymath}
	{\mathcal P}(V_0)/{\mathcal I}\equiv V_\Pi.
\end{displaymath} 

The finite-dimensionality assumption on $V_\Pi$ then implies that
${\mathcal I}$ contains all homogeneous polynomials of a degree
greater than or equal to some fixed degree, say $d_0$. Since there are
only a finite number of linearly independent homogeneous polynomials
in  ${\mathcal P}(V_0)$ of degree $d_0$, it follows that there exists
a finite number of elements $p_1,p_2,\dots,p_j$ in  ${\mathcal P}(V_0)$ (these may be chosen for
instance as  highest weight vectors) such that if ${\mathcal I}\langle
p_1,p_2,\dots, p_j\rangle$ denotes the  ${\mathcal S}( {\mathfrak
  p}^-)\times_s {\mathfrak g}_0$ submodule generated by these elements, then

\begin{displaymath}
  V_\Pi\equiv  {\mathcal P}(V_0)/{\mathcal I}\langle
p_1,p_2,\dots, p_j\rangle.
\end{displaymath}

We assume that the number $j$ of polynomials is minimal. 

\medskip

Once the elements $p_1,p_2,\dots,p_j$ are known, it is possible to
construct the whole graph as above. In particular,  {\bf the sinks in
  $V_\Pi$ can be directly determined from this}. See
Proposition~\ref{clearprop} below for a simple example that indicates
how. In case $\dim V_0$ is large the task, of course, will be more
cumbersome. 

\medskip

\begin{Ex} As is well known, we have that $${\mathfrak p}^-\otimes \textrm{spin}(\frac12,0)=\textrm{spin}(1,\frac12)\oplus \textrm{spin}(0,\frac12).$$
If we mod out by all second order polynomials, and possibly one of the
first order polynomial representations, we get the following 3
indecomposable representations: 
\begin{itemize}
	\item $\textrm{spin}(\frac12,0)\rightarrow
          \textrm{spin}(0,\frac12)$. This 4-dimensional representation
          comes from the the defining representation.
	\item $\textrm{spin}(\frac12,0)\rightarrow
          \textrm{spin}(1,\frac12)$. This is an 8-dimensional
          representation. 
	\item $\textrm{spin}(\frac12,0)\rightarrow
          \textrm{spin}(0,\frac12), \textrm{spin}(1,\frac12)$. This is
          a 10-dimensional representation which includes the two
          former. 
\end{itemize}
\medskip

Proceeding analogously, $${\mathfrak p}^-\otimes \textrm{spin}(1,0)=\textrm{spin}(\frac32,\frac12)\oplus \textrm{spin}(\frac12,\frac12)$$
leads to inequivalent representations in dimensions 7,11,15. 

Similarly,
$${\mathfrak p}^-\otimes
\textrm{spin}(\frac12,\frac12)=\textrm{spin}(0,0)\oplus
\textrm{spin}(1,0)\oplus \textrm{spin}(0,1)\oplus \textrm{spin}(1,1)$$

leads to indecomposable representations in dimensions 5, 7, 8, 10, 11,
13, 16, 17, 19, 20. Several dimensions here carry a number of
inequivalent representations. 

Together with duals of these or versions obtained as mirror images by
interchanging the spins, these exhaust all the known representations
in dimensions less than or equal to 8 with the exception of a
6-dimensional representation which we describe in Example~\ref{ex2}. 
\end{Ex}

\medskip

\subsection{Special case: Ideals in ${\mathcal U}({\mathfrak p}^-)$} It is well-known that there is a decomposition
\begin{displaymath}
	{\mathcal U}({\mathfrak p}^-)=\oplus W_{r,s}
\end{displaymath}
into ${\mathfrak g}_0$ representations, where the subspace $W_{r,s}$
may be defined through its highest weight vector, say $z_1^r\det
z^s$. This is possible since each representation occurs with
multiplicity one. Denote this representation simply by $[r,s]$. It is
elementary to see that the action of ${\mathfrak p}^-$ on the left on
${\mathcal U}({\mathfrak p}^-)$, when expressed in terms of
representations, is given as follows. All arrows represent non-trivial
maps. 

\begin{equation}
	\begin{array}{cccccccc}\\\downarrow\\\left[1,0\right]\\\downarrow&\searrow\\\left[2,0\right]&&\left[0,1\right]\\\downarrow&\searrow&\downarrow\\\left[3,0\right]&&\left[1,1\right]\\\downarrow&\searrow&\downarrow&\searrow\\\left[4,0\right]&&\left[2,1\right]&&\left[0,2\right]\\\downarrow&\searrow&\downarrow&\searrow&\downarrow\\\left[5,0\right]&&\left[3,1\right]&&\left[1,2\right]\\\downarrow&\searrow&\downarrow&\searrow&\downarrow&\searrow\\\left[6,0\right]&&\left[4,1\right]&&\left[2,2\right]&&\left[0,3\right]\\\vdots&&\vdots&&\vdots&&\vdots 
\end{array}
	\label{tree}
\end{equation}

Any ideal ${\mathcal I}\subseteq {\mathcal U}({\mathfrak p}^-)$ that
has finite codimension must clearly contain some $z_1^{r_1}$ for some
minimal $r_1\in {\mathbb N}$ (we omit the trivial case of codimension
0). Since we are assuming that the ideals are ${\mathfrak g}_0$
invariant, if some other element $z_1^{r_2}p(\det z)$ is in the ideal
then first of all we can assume $r_2<r_1$ and secondly we can assume
that the polynomial $p$ is homogeneous; $p(\det z)=\det z^{s_2}$ for
some $s_2>0$. The latter inequality follows by the minimality of
$r_1$. 

Thus the following is clear:
\begin{Prop} \label{clearprop}
Any ${\mathfrak g}_0$ invariant ideal ${\mathcal I}\subseteq {\mathcal
  U}({\mathfrak p}^-)$ of finite codimension is of the form 
\begin{equation}
	{\mathcal I}={\mathfrak g}_0\cdot \langle z_1^{r_1}\det
        z^{s_1},z_1^{r_2}\det z^{s_2},\cdots, z_1^{r_t}\det
        z^{s_t}\rangle 
	\label{eqw}
\end{equation}
for some positive integers $r_1>r_2>\dots>r_t$ and integers
$0=s_1<s_2<\cdots<s_t$. If the set is minimal, then  furthermore 
$$\forall j=2,3,\dots, t: s_1+s_2+\dots + s_j\leq r_1-r_j.$$
Any set of such integers determine an invariant ideal of finite
codimension. The sinks in the quotient module are
$$z_1^{r_1-s_2}\det z^{s_2-1},z_1^{r_1-s_2-s_3}\det z^{s_3-1}\dots
z_1^{r_1-s_2-\dots-s_t}\det z^{s_t-1},\textrm{ and } \det z^{s_t+r_t-1}.$$
\end{Prop}

\medskip
\begin{Ex}\label{ex2} If we mod out by all second order polynomials,
  that is by $z_1^2$ and $\det z$, we get the 5-dimensional
  indecomposable representation $$\textrm{spin}(0,0)\rightarrow
  \textrm{spin}(\frac12,\frac12).$$ If we instead mod out by $z_1^2$
  and $z_1\det z$ we get the 6-dimensional indecomposable
  representation   
$$\textrm{spin}(0,0)\rightarrow \textrm{spin}(\frac12,\frac12)\rightarrow \textrm{spin}(0,0).$$
\end{Ex}

\medskip

\begin{Ex} The representations determined by ideals are easily written
  down, though some finer details from the representation theory of
  $su(2)$ will have to be invoked to get the precise form. In simple
  examples like the last in the previous example, everything follows
  immediately since there is no need to be precise about the relative
  scales in the 3 spaces: 

\begin{eqnarray}
{\mathfrak p}^-\ni \underline{p}=(p_1,p_2,p_3,p_4)\mapsto
\left(\begin{array}{cccccc}0&0&0&0&0&0\\p_1&0&0&0&0&0\\p_2&0&0&0&0&0\\p_3&0&0&0&0&0\\p_4&0&0&0&0&0\\0&p_4&-p_3&-p_2&p_1&0    
\end{array}\right).
\end{eqnarray}
An element $(u,v)$ in the diagonal subgroup $G_0$ (see
Section~\ref{fd}) acts as $0\oplus u\otimes (v^t)^{-1}\oplus 0$. 

Notice that the matrix with just $p_1$ corresponds to a map which
sends the constant 1 to the polynomial $p_1z_1$ and sends the
polynomial $z_4$ to $p_1\det z$ and all other first order polynomials
$z_1,z_2,z_3$ to 0. 

\end{Ex}

\medskip

\subsection{Two sources and 2 sinks}

We here consider the Poincar\'e algebra.

Consider the situation previously depicted under `Simple situations'
where there is one source and two sinks. The resulting representations
may be written as

\begin{equation}\label{repg1}
  \left(\begin{array}{cc}0&0\\w&0\end{array}\right)\mapsto
  \left(\begin{array}{ccc}0&0&0\\F(w)&0&0\\G(w)&0&0
\end{array}\right)\; ,
  \left(\begin{array}{cc}u&0\\0&v\end{array}\right)\mapsto
  \left(\begin{array}{ccc}\tau_1(u,v)&0&0\\0&\tau_2(u,v)&0\\0&0&\tau_3(u,v) 
\end{array}\right). 
\end{equation}

With this one can easily write down a representation with 2 sources
and 2 sinks, indeed a 4-parameter family given by elements
$(\alpha,\beta,\gamma,\delta)\in {\mathbb C}^4$:

\begin{equation}\label{repg2}  
  \left(\begin{array}{cccc}0&0&0&0\\0&0&0&0\\\alpha\cdot
      F(w)&\beta\cdot F(w)&0&0\\\gamma \cdot G(w)&\delta\cdot
      G(w)&0&0 
\end{array}\right)\textrm{ resp. }
\left(\begin{array}{cccc}\tau_1(u,v)&0&0&0\\0&\tau_1(u,v)&0&0\\0&0&\tau_2(u,v)&0\\0&0&0&\tau_3(u,v)
\end{array}\right). 
\end{equation}

\medskip

In this case, there is a continuum of inequivalent representations and they are
generically indecomposable. This lead to a continuum already
in dimension 12 where the two sources are equal and 2-dimensional -
say spin($\frac12$,0), and the two sinks are spin($1$,$\frac12$) and
spin(0,$\frac12$) 
or, also in dimension 12, the 2 sources are the 4-dimensional
spin($\frac12$,$\frac12$), and the sinks are spin($1$,0) and
spin(0,0), or in  
dimension 16 where one is the 2-dimensional spin ($\frac12$,0) and the 
other is the 6-dimensional spin($\frac12$,$1$) and the targets are
spin(0,$\frac12$) and spin(1,$\frac12$). The moduli space in these
cases is ${\mathbb C}{\mathbb P}^1$. Specifically, the indecomposable
is determined by a point $(p,q) \in {\mathbb C}{\mathbb P}^1\times
{\mathbb C}{\mathbb P}^1$ giving relative scales on the arrows. Here,
$p\equiv(\alpha,\beta)$ and $q\equiv(\gamma,\delta)$ in the above
representation. Two such points $(p_1,q_1)$ and $(p_2,q_2)$, are
equivalent if there is an element $g\in GL(2,{\mathbb C})$ such that
$(p_2,q_2)=(gp_1,gp_2)$.

\bigskip

\section{Supersymmetry}

The previous considerations are now extended to the supersymmetric
setting as follows: Let $H_1, H_2$, and $H_3$ be reductive matrix Lie
groups with Lie algebras ${\mathfrak h}_1, {\mathfrak h}_2$, and
${\mathfrak h}_3$, respectively. We assume that possible abelian
ideals are represented by semi-simple elements. We consider an
irreducible representation of each of these Lie algebras;  $V_1,V_2$,
and $V_3$. We identify the representation with the space in which it
acts. We denote furthermore the dual representation of a
representation $V$ by $V^\prime$ (this is the ${\mathbb C}$ linear
dual). Let  
\begin{eqnarray}W_1= \hom(V_1,V_2)\equiv V_1^\prime\otimes V_2 &;&W_2=
  \hom(V_2,V_3)\equiv V_2^\prime\otimes V_3\\\qquad \textrm{and}
  \qquad Z= \hom(V_1,V_3)\equiv V_1^\prime\otimes V_3. 
\end{eqnarray}    

The Lie superalgebra ${\mathfrak g}_{super}={\mathfrak
  g}_{super}({\mathfrak h}_1, {\mathfrak h}_2,{\mathfrak
  h}_3,V_1,V_2,V_3)$  is defined as 
\begin{eqnarray}
  &{\mathfrak g}_{super}=\\&\left\{\left(\begin{array}{ccc}a&0&0\\w_1&g&0\\z&w_2&b\end{array}\right)
\mid 
    a\in{\mathfrak h}_1\; , g \in
    {\mathfrak h}_2\; , b\in {\mathfrak h}_3\; , w_1\in W_1\;,w_2\in W_2\;, \textrm{ and }z\in Z  \right\}.\nonumber
\end{eqnarray}The odd part is given as 

\begin{eqnarray}
  {\mathfrak g}_{super}^1=\left\{\left(\begin{array}{ccc}0&0&0\\w_1&0&0\\0&w_2&0\end{array}\right)
\mid w_1\in W_1\;,\textrm{ and }w_2\in W_2\;\right\}.\nonumber
\end{eqnarray}

Let  
\begin{eqnarray}
  {\mathfrak n}_{super}=\left\{\left(\begin{array}{ccc}0&0&0\\w_1&0&0\\z&w_2&0\end{array}\right)
\mid w_1\in W_1\;,w_2\in W_2\; ,\textrm{ and }z\in Z\right\}.\nonumber
\end{eqnarray}
and
\begin{eqnarray}
  {\mathfrak g}^r_{super}=\left\{\left(\begin{array}{ccc}a&0&0\\0&g&0\\0&0&b\end{array}\right)
\mid
    a\in{\mathfrak h}_1\; , g \in
    {\mathfrak h}_2\; , \textrm{ and }b\in {\mathfrak h}_3\;    \right\}.\nonumber
\end{eqnarray}

Obviously, ${\mathfrak n}_{super}$ is a maximal nilpotent ideal and
${\mathfrak g}_{super}^r$ is the reductive part.  

We let \begin{eqnarray}
  {\mathfrak p}^-=\left\{\left(\begin{array}{ccc}0&0&0\\0&0&0\\z&0&0\end{array}\right)
\mid z\in Z\right\}.\nonumber
\end{eqnarray}

Then ${\mathfrak g}_{super}^0={\mathfrak g}^r_{super}\oplus {\mathfrak p}^-$.

We then have the following super algebraic generalization of Theorem~\ref{lem1}:

\begin{Thm}Consider a finite-dimensional representation $V_{super}$ of
  ${\mathfrak g}_{super}$. Then there is a flag of subspaces
  $\{0\}=W_0\subsetneq W_1\subsetneq\cdots\subsetneq W_{r-1}\subsetneq
  W_r=V_{super}$ such that each $W_i$ is invariant and completely
  reducible under ${\mathfrak g}_{super}^0$ and such that ${\mathfrak
    n}_{super}$ maps  $W_i$ into $W_{i-1}$ for $i=1,\dots,r$. 
\end{Thm}

Thus, the previous treatment with directed graphs and ideals carry
over. Naturally, the picture gets more complicated with the odd
generators.

The most simple thing to consider would be the ${\mathfrak g}_{super}$
module ${\mathcal U}({\mathfrak n}_{super})$, but even here the
situation is complex, though in principle tractable. 

Consider as an example the case of the simplest super Poincar\'e algebra,

\begin{eqnarray} &{\mathfrak
    g}_{super}^{P}=\\&\left\{\left(\begin{array}{ccc}a&0&0\\w_1&0&0\\z&w_2&b\end{array}\right) 
\mid
    a,b\in sl(2,{\mathbb C}) , w_1 \in M(1,2), w_2\in M(2,1)\;, \textrm{ and }z\in M(2,2)  \right\}.\nonumber
\end{eqnarray} 

\medskip

\medskip

Let $W_1=M(1,2)$ and $W_2=M(2,1)$. We  number the spaces $$\begin{array}{ccc}1&W_2&W_2\wedge
  W_2\\W_1&W_1\wedge W_2&W_1\wedge W_2\wedge W_2\\W_1\wedge
  W_1&W_1\wedge W_1\wedge W_2&W_1\wedge W_1\wedge W_2\wedge W_2\end{array}=\begin{array}{ccc}1&2&3\\4&5&6\\7&8&9\end{array}.$$ 

We then have that \begin{equation}
	{\mathcal U}({\mathfrak n}_{super})=\sum_{i=1}^9 {\mathcal
          U}({\mathfrak p}^-)\otimes {\mathcal U}({\mathfrak
          q}_{super}^1)_i. 
	\label{nsu}
\end{equation}

\medskip

Each of the 9 summands is invariant under ${\mathfrak g}^r_{super}$.
The representations corresponding to this are given right below
. Here, $n,d$ are independent non-negative integers (in a few obvious
cases, n must furthermore be non-zero). The labels $\uparrow$ and
$\downarrow$ may be taken just as part of a short hand notation that
are defined by the stated equations. They are  listed  here for
convenience even though they are not used directly. One can ascertain
useful information from them about how the various pieces occur in the
tensor products of $su(2)$ representations.  

\begin{eqnarray*}
1\left[n,n,d\right]&\oplus&\\
2\uparrow\left[n,d\right]=2\left[n-1,n,d\right]&\oplus&2\downarrow\left[n,d\right]=2\left[n+1,n,d\right]\\
3\left[n,n,d\right]&\oplus&\\
4\uparrow\left[n,d\right]=4\left[n,n-1,d\right]&\oplus&4\downarrow\left[n,d\right]=4\left[n,n+1,d\right]\\
5\uparrow\uparrow\left[n,d\right]=5\left[n-1,n-1,d\right]&\oplus&5\uparrow\downarrow\left[n,d\right]=5\left[n-1,n+1,d\right]\\
5\downarrow\uparrow\left[n,d\right]=5\left[n+1,n-1,d\right]&\oplus&5\downarrow\downarrow\left[n,d\right]= 5\left[n+1,n+1,d\right]\\
6\uparrow\left[n,d\right]=6\left[n,n-1,d\right]&\oplus&6\downarrow\left[n,d\right]=6\left[n+1,n,d\right]\\
7\left[n,n,d\right]&\oplus&\\
8\uparrow\left[n,d\right]=6\left[n-1,n,d\right]&\oplus&8\downarrow\left[n,d\right]=8\left[n+1,n,d\right]\\
&\oplus&9\left[n,n,d\right] 
\end{eqnarray*}

A further complication is that there are representations in different
spaces that are equivalent under ${\mathfrak g}^r_{super}$:

\begin{eqnarray*}
	\label{relations}
8\uparrow\left[n,d\right]&\leftrightarrow&4\uparrow\left[n-1,d+1\right]\\
8\downarrow\left[n,d\right]&\leftrightarrow&4\left[n+1,d\right]\\
5\downarrow\downarrow\left[n,d\right]&\leftrightarrow&1\left[n+1,d\right]\\
5\uparrow\uparrow\left[n,d\right] &\leftrightarrow&1\left[n-1,d+1\right]\\
5\downarrow\downarrow\left[n-1,d+1\right],5\uparrow\uparrow\left[n+1,d\right]&\leftrightarrow&9\left[n,d\right]\\
5\downarrow\downarrow\left[n-1,d+1\right],5\downarrow\downarrow\left[n+1,d\right]&\leftrightarrow&9\left[n,d\right]\\
6\uparrow\left[n,d\right]&\leftrightarrow&2\downarrow\left[n-1,d+1\right]\\
6\uparrow\left[n,d\right]&\leftrightarrow&2\downarrow\left[n-1,d+1\right]\\
6\downarrow\left[n,d\right]&\leftrightarrow&2\uparrow\left[n+1,d\right]\\
6\downarrow\left[n,d\right]&\leftrightarrow&2\uparrow\left[n+1,d\right]
\end{eqnarray*}

To each finite-dimensional representation $V_r$ of  ${\mathfrak
  g}^r_{super}$ (may be reducible), the general object of interest is
the left module 
\begin{equation}
	{\mathcal U}({\mathfrak n}_{super})\cdot V_r=\sum_{i=1}^n
        {\mathcal U}({\mathfrak p}^-)\otimes {\mathcal U}({\mathfrak
          q}_{super}^1)_i\cdot V_r. 
	\label{nsu1}
\end{equation}

To further analyze this we have to choose a PBW-type basis. We will do this in the indicated fashion with ${\mathcal U}({\mathfrak p}^-)$ to the left and with furthermore $W_1$ always to the left of $W_2$.

\begin{Ex}Assume that ${\mathcal U}({\mathfrak p}^-)$ acts trivially
  on the space $V_r$. The resulting module is then  
\begin{equation}\bigwedge({\mathfrak g}^1_{super})\cdot V_r,
	\label{triv}
\end{equation}
or some of the subrepresentations thereof. The exterior algebra
$\bigwedge({\mathfrak g}^1_{super})$ occurs because the $W_1$, $W_2$
anticommute in the considered quotient. 
\end{Ex}

\medskip

Observe that in the sum (\ref{nsu1}) the summand
\begin{equation}
	{\mathcal U}_{2,3,5,6,8,9}({\mathfrak n}_{super})\cdot
        V_r=\sum_{i=2,3,5,6,8,9} {\mathcal U}({\mathfrak p}^-)\otimes
        {\mathcal U}({\mathfrak q}_{super}^1)_i\cdot V_r 
	\label{nsu2}
\end{equation}
is invariant. We may then pass to a general subclass of indecomposable
modules by first taking the quotient by this. The vector space that
results is   
\begin{equation}
	{\mathcal U}_{rest}({\mathfrak n}_{super})\cdot V_r=\sum_{i=1,4,7} {\mathcal U}({\mathfrak p}^-)\otimes {\mathcal U}({\mathfrak q}_{super}^1)_i\cdot V_r.
	\label{nsu3}
\end{equation}

If we let ${\mathcal U}({\mathfrak p}^-)_{\geq s}$ be the ideal generated
by all homogeneous elements of degree $s$, it is easy to see that for
each $s=0,1,2,\dots$, the space  

\begin{equation}
	{\mathcal U}_{rest}^s({\mathfrak n}_{super})\cdot V_r=
        {\mathcal U}({\mathfrak p}^-)_{s+2}\cdot  V_r \oplus {\mathcal
          U}({\mathfrak p}^-)_{s+1}\cdot W_1\cdot V_r\oplus {\mathcal
          U}({\mathfrak p}^-)_{s}\cdot (W_1\wedge W_1)\cdot V_r 
	\label{nsu4}
\end{equation}
is invariant.

\begin{Ex} Let $V_r^1$ be the irreducible 2-dimensional representation
  which is only non-trivial on the ${\mathfrak h}_1$ piece of the
  reductive part.  The defining representation of ${\mathfrak
    g}_{super}^{P}$ is a subrepresentation of the quotient $${\mathcal
    U}_{rest}({\mathfrak n}_{super})\cdot V^1_r/{\mathcal
    U}_{rest}^0({\mathfrak n}_{super})\cdot V^1_r.$$ Indeed, we just
  have to limit ourselves further by removing two appropriate
  ${\mathfrak g}_{super}^r$ representations. 
\end{Ex}

Returning to the more general situation, let us assume from now on
that $V_r$ is the trivial 1-dimensional module. We are thus left with
the space  

\begin{equation}
	{\mathcal U}_{rest}({\mathfrak n}_{super})= {\mathcal
          U}({\mathfrak p}^-) \oplus {\mathcal U}({\mathfrak
          p}^-)\cdot W_1\oplus {\mathcal U}({\mathfrak p}^-)\cdot
        (W_1\wedge W_1). 
	\label{nsu5}
\end{equation}

The general form of the representation is (we give only the $W_1,W_2$ operators)
\begin{equation}
	\left(\begin{array}{cccc}0&(w_2^1p_{11}+w_2^2p_{21})1_{\mathcal P}&(w_2^1p_{12}+w_2^2p_{22})1_{\mathcal P}&0\\w_1^11_{\mathcal P}&0&0&-(w_2^1p_{12}+w_2^2p_{22})1_{\mathcal P}\\w_1^21_{\mathcal P}&0&0&(w_2^1p_{11}+w_2^2p_{21})1_{\mathcal P}\\0&-w_1^21_{\mathcal P}&w_1^11_{\mathcal P}&0 
\end{array}\right) .
	\label{genfo}
\end{equation}
Here, each block corresponds to a space of polynomials. The operators
$p_{i,j}$ are multiplication operators in ${\mathcal U}({\mathfrak
  p}^-)$ - and hence also such operators in the given space. Notice
that they increase the degree of the target by 1. The symbol
$1_{\mathcal P}$ denotes the identity operator. 

\medskip

The finer details are given as follows, where the arrows point upwards
come from $W_2$ and those pointing downwards come from $W_1$.

\begin{equation}
\begin{array}{cccccccccccccccc}&&1\left[n,d+1\right]&&\\&\nearrow&&\nwarrow&\\4\downarrow\left[n-1,d+1\right]&&&&4\uparrow\left[n+1,d\right]\\&\nwarrow&&\nearrow&\\
&&7\left[n,d\right]&&
\end{array}
\end{equation}

\begin{equation}
\begin{array}{cccccccccccccccc}&&1\left[n,d\right]&\\&\swarrow&&\searrow\\4\downarrow\left[n,d\right]&&&&4\uparrow\left[n,d\right]
  \\&\nwarrow&&\nearrow\\
&&7\left[n,d\right] 
\end{array}
\end{equation}

\medskip

Any invariant ideal ${\mathcal I}_{super}$ in ${\mathcal
  U}_{rest}({\mathfrak n}_{super})$ contains a sum of the form  

\begin{equation}
	 {\mathcal I}_1({\mathfrak p}^-) \oplus {\mathcal
           I}_4({\mathfrak p}^-)\cdot W_1\oplus {\mathcal
           I}_7({\mathfrak p}^-)\cdot (W_1\wedge W_1), 
	\label{nsu6}
\end{equation}
where ${\mathcal I}_1({\mathfrak p}^-)\subseteq {\mathcal
  I}_4({\mathfrak p}^-)\subseteq {\mathcal I}_7({\mathfrak
  p}^-)\subseteq {\mathcal U}({\mathfrak p}^-)$ are ${\mathfrak p}^-$
ideals. These are precisely the ideals determined in Section~1. 

\begin{Prop}
If we have a representation $7\left[n,d\right]\in {\mathcal
  I}_7({\mathfrak p}^-)$ then $4\left[n,d+1\right]\in {\mathcal
  I}_4({\mathfrak p}^-)$ and $1\left[n,d+1\right]\in {\mathcal
  I}_1({\mathfrak p}^-)$. Furthermore, $4\uparrow\left[n,d\right],
4\downarrow\left[n,d\right]\in {\mathcal I}_{super}$.  In particular,
$${\mathfrak p}^-\cdot {\mathfrak p}^-\cdot{\mathcal I}_7({\mathfrak
  p}^-)\subseteq {\mathcal I}_1.$$
\end{Prop}

This result in principle solves the problem but there are still
extremely many cases - even if we start with an ideal ${\mathcal I}_7$
and ask for how many configurations of the ideals ${\mathcal I}_1,
{\mathcal I}_4$ that are possible. We refrain from pursuing this
further and just give a low-dimensional example. 

\begin{Ex}Let ${\mathcal I}_7={\mathcal I}\langle z_1\rangle$. The
  following list is exhaustive and each case occurs. 
\begin{itemize}
	\item ${\mathcal I}_1={\mathcal I}\langle z_1\rangle$ then
          ${\mathcal I}_4={\mathcal I}\langle z_1\rangle$. 
	\item ${\mathcal I}_1={\mathcal I}\langle z^2_1\rangle$ then
          either ${\mathcal I}_4={\mathcal I}\langle z^2_1\rangle$,
          ${\mathcal I}_4={\mathcal I}\langle z^2_1, \det z\rangle$,
          or ${\mathcal I}_4={\mathcal I}_7$.  
	\item ${\mathcal I}_1={\mathcal I}\langle z^2_1, \det
          z\rangle$ then ${\mathcal I}_4={\mathcal I}\langle z^2_1,
          \det z\rangle$. 
	\item ${\mathcal I}_1={\mathcal I}\langle z^3_1, \det
          z\rangle$ then ${\mathcal I}_4={\mathcal I}\langle z^2_1,
          \det z\rangle$. 
	\item ${\mathcal I}_1={\mathcal I}\langle z^3_1, z_1\det
          z\rangle$ then ${\mathcal I}_4={\mathcal I}\langle z^2_1,
          \det z\rangle$ or ${\mathcal I}_4={\mathcal I}^2_1$.   
\end{itemize} 
\end{Ex}

\end{document}